\RequirePackage{ifpdf}
\ifpdf % We are running pdfTeX in pdf mode
\documentclass[pdftex]{sigma}
\else
\documentclass{sigma}
\fi

\usepackage{mathrsfs}

\newcommand{\re}{\mathbb{R}}

\newcommand{\be}{{}_0\mathscr{F}_1}

\begin{document}

\renewcommand{\thefootnote}{$\star$}

\renewcommand{\PaperNumber}{075}

\FirstPageHeading

\ShortArticleName{Generalized Bessel function of Type $D$}

\ArticleName{Generalized Bessel function of Type $\boldsymbol{D}$\footnote{This paper is a contribution to the Special
Issue on Dunkl Operators and Related Topics. The full collection
is available at
\href{http://www.emis.de/journals/SIGMA/Dunkl_operators.html}{http://www.emis.de/journals/SIGMA/Dunkl\_{}operators.html}}}

\Author{Nizar DEMNI}

\AuthorNameForHeading{N. Demni}

\Address{SFB 701, Fakult\"at f\"ur Mathematik, Universit\"at Bielefeld, Deutschland}
\Email{\href{mailto:demni@math.uni-bielefeld.de}{demni@math.uni-bielefeld.de}}

\ArticleDates{Received July 01, 2008, in f\/inal form October 24,
2008; Published online November 04, 2008}

\Abstract{We write down the generalized Bessel function associated with the root system of type $D$ by means of multivariate hypergeometric series. Our hint comes from the particular case of the Brownian motion in the Weyl chamber of type $D$.}


\Classification{33C20; 33C52; 60J60; 60J65}

\renewcommand{\thefootnote}{\arabic{footnote}}
\setcounter{footnote}{0}

\section{Root systems and related processes}

We refer the reader to \cite{Hum} for facts on root systems. Let $(V,\langle\cdot\rangle)$ be an Euclidean space of f\/inite dimension $m \geq 1$. A \emph{reduced} root system $R$ is a f\/inite set of non zero vectors in $V$ such that
\begin{itemize}\itemsep=0pt
\item[1)] $R \cap \re \alpha = \{\alpha,-\alpha\}$ for all $\alpha \in R$,
\item[2)] $\sigma_{\alpha}(R) = R$,
\end{itemize} where $\sigma_{\alpha}$ is the ref\/lection with respect to the hyperplane $H_{\alpha}$ orthogonal to $\alpha$
\begin{gather*}
\sigma_{\alpha} (x) =  x - 2\frac{\langle \alpha, x\rangle}{\langle\alpha,\alpha\rangle} \alpha , \qquad x \in V.
 \end{gather*}
A simple system $S$ is  a basis of $\textrm{span}(R)$ which induces a total ordering in $R$. A  root $\alpha$ is positive if it is a positive linear combination of elements of $S$. The set of positive roots is called a positive subsystem and is denoted by $R_+$. The (f\/inite) ref\/lection group $W$ is the group generated by all the ref\/lections $\sigma_{\alpha}$ for $\alpha \in R$.
Given a root system $R$  with positive and simple systems $R_+$, $S$,  def\/ine the {\it positive Weyl chamber} $C$ by:
\begin{gather*}
C := \{x \in V, \, \langle \alpha , x \rangle > 0 \  \forall \, \alpha \in R_+\} = \{x \in V, \  \langle \alpha , x \rangle > 0 \  \forall \, \alpha \in S\}
 \end{gather*} and $\overline{C}$, $\partial C$ its closure and boundary respectively. One of the most important properties is that the convex cone $\overline{C}$ is a fundamental domain, that is, each $\lambda \in V$ is conjugate to one and only one $\mu \in \overline{C}$.
Processes related to root systems have been of great interest during the last decade and notably the so-called Dunkl process~\cite{Ros0}. This $V$-valued process was deeply studied in a~sequence of papers by Gallardo and Yor \cite{Gal1,Gal2,Gal3,GG1} and Chybiryakov~\cite{Chy} and its projection on the Weyl chamber gives rise to a dif\/fusion known as the $W$-invariant or radial Dunkl process. The generator of the latter process acts on  $C_c^2(\overline{C})$ as
\begin{gather*}
\mathscr{L}_ku(x) = \frac{1}{2}\Delta u(x) + \sum_{\alpha \in R_+}k(\alpha)\frac{\langle \alpha,\nabla u(x)\rangle}{\langle \alpha, x\rangle}
\end{gather*}
with $\langle \alpha,\nabla u(x)\rangle = 0$ whenever $\langle \alpha,x\rangle = 0$, where $k$ is a positive function def\/ined on the orbits of the $W$-action on $V$ and is constant on each orbit. Such a function is called a multiplicity function. The semi group density (with respect to the Lebesgue measure on $V$) of the radial Dunkl process is written as~\cite{Ros0}
\begin{gather}
\label{sg}
p_t^k(x,y) = \frac{1}{c_kt^{\gamma + m/2}}e^{-(|x|^2+|y|^2)/2t}D_k^W\left(\frac{x}{\sqrt t},\frac{y}{\sqrt t}\right)\omega_k(y)^2
\end{gather}
for $x,y \in C$, where
\begin{gather*}
\gamma = \sum_{\alpha\in R_+}k(\alpha),\qquad \omega_k(y) = \prod_{\alpha \in R_+}\langle \alpha,y \rangle^{k(\alpha)}.
\end{gather*}
The kernel $D_k^W$ is def\/ined by
\begin{gather*}
D_k^W(x,y) = \sum_{w \in W}D_k(x,wy)
\end{gather*}
where $D_k$ is the non symmetric Dunkl kernel~\cite{Ros1}, and is up to the constant factor $1/|W|$ the so-called \emph{generalized Bessel function} since it reduces to the (normalized) Bessel function in the rank-one case~$B_1$~\cite{Ros1}. It was identif\/ied in \cite{Bak} with a multivariate hypergeometric series for both root systems of types~$A$ and~$B$. Since the $C$-type root system is nothing but the dual root system of $B$, it remains to write down this kernel for the root system of type $D$ which is the remaining inf\/inite family of irreducible root systems associated with Weyl groups. This was behind our motivation to write the present paper and our main result is stated as

\begin{theorem}\label{theorem1}
For the $D$-type root system, the generalized Bessel function writes
\begin{gather*}
\frac{1}{|W|}\sum_{w \in W} D_k(x,wy) = \prod_{i=1}^m \left(\frac{x_iy_i}{2}\right){}_0F_1^{(1/k_1)}\left(q+\frac{1}{2}, \frac{x^2}{2},\frac{y^2}{2}\right)+
{}_0F_1^{(1/k_1)}\left(q-\frac{1}{2}, \frac{x^2}{2},\frac{y^2}{2}\right),
\end{gather*}
where  ${}_0F_1^{(1/k_1)}$ is the multivariate hypergeometric series of Jack parameter $1/k_1$, $x := (x_1^2,\dots,$ $x_m^2)$, $k = k_1 > 0$ is the value of the multiplicity function and $q := 1+(m-1)k_1$.\end{theorem}

Our strategy is easily explained as follows: according to Grabiner~\cite[p.~21]{Gra}, the Brownian motion in the Weyl Chamber is a radial Dunkl process of multiplicity function $k \equiv 1$, that is $k(\alpha)= 1$ for all $\alpha \in R$. Indeed, it is shown in \cite{Scha} that the radial Dunkl process, say $X^W$, is the unique strong solution of
\begin{gather*}
dY_t = dB_t + \sum_{\alpha \in R_+}\frac{k(\alpha)dt}{\langle\alpha,Y_t\rangle}, \qquad Y_0 \in \overline{C},
\end{gather*}
where $B$ is a $m$-dimensional Brownian motion and $k(\alpha) >0$ for all $\alpha \in R$.  Hence, by Grabiner's result on Brownian motions started at $C$ and killed when they f\/irst reach $\partial C$ together with the Doob's transform~\cite{Rev} applied with the harmonic function~\cite[p.~19]{Gra}
\begin{gather*}
h(y) := \prod_{\alpha \in R_+}\langle \alpha,y\rangle,
\end{gather*}
the semi group density the Brownian motion in the Weyl chamber is given by
\begin{gather}\label{BmWC}
p_t^1(x,y) = \frac{h(y)}{h(x)}\sum_{w \in W}\det(w) N_{t,m}(wy-x),
\end{gather}
where $N_{t,m}$ is the $m$-dimensional heat kernel given by
\begin{gather*}
N_{t,m}(x) = \frac{1}{(2\pi t)^{m/2}}e^{-|x|^2/2t}, \qquad |x|^2 = \langle x,x\rangle, \qquad x \in V.
\end{gather*}
On the one hand and for the irreducible root systems $A$, $B$, $D$, the sum over $W$ in (\ref{BmWC}) was expressed in~\cite{Gra} as determinants. On the other hand, the generalized Bessel function was expressed for $A$ and $B$-types root systems via multivariate hypergeometric series of two arguments and of Jack parameter which equals the inverse of one of the multiplicity values (see the end of~\cite{Bak}).  When this parameter equals one, the multivariate series takes a determinantal form~\cite{Gro} which agrees with Grabiner's result. As a matter of fact, we will start from the expression obtained in \cite{Gra} for the Brownian motion in the Weyl chamber of type $D$ then express it as a multivariate series of a Jack parameter~$1$. Once we did, we will generalize the result to an arbitrary positive Jack parameter using an analytical tool known as the shift principle. Before proceeding, we shall recall some facts on multivariate hypergeometric series and investigate the cases of $A$ and $B$-types root systems.

\section{Mutlivariate series and determinantal representations}

We refer the reader to \cite{Bak,Bee,Kan} and references therein for facts on Jack polynomials and multi\-variate hypergeometric series. Let $\tau$ be a partition of length $m$, that is a sequence of positive integers $\tau_1 > \cdots > \tau_m$. Let $\alpha > 0$, then the Jack polynomial $C_{\tau}^{(\alpha)}$ of Jack parameter $\alpha$ is the unique (up to a normalization) symmetric homogeneous eigenfunction of the operator
\begin{gather*}
\sum_{i=1}^mx_i^2\partial_i^2  + \frac{2}{\alpha}\sum_{i=1}^m\frac{x_i^2}{x_i-x_j}\partial_i
\end{gather*}
corresponding to the eigenvalue
\begin{gather*}
\rho_{\tau} = \sum_{i=1}^m \tau_i[\tau_i - (2/\alpha)(i-1)] + |\tau|(m-1), \qquad |\tau| = \tau_1 +\cdots + \tau_m.
\end{gather*}
The normalization we adopt here is
\begin{gather*}
(x_1+\cdots +x_m)^n = \sum_{\tau}C_{\tau}^{(\alpha)}(x),\qquad n \geq 0,
\end{gather*}
where the sum is taken over all the partitions of weight $|\tau| = n$ and length $m$. For $\alpha = 2$ it reduces to the zonal polynomial~\cite{Muir} while for $\alpha = 1$ it f\/its the Schur polynomial~\cite{Mac}. Let $p,q \in \mathbb{N}$, then the multivariate hypergeometric series of two arguments is def\/ined by
\begin{gather*}
{}_pF_q^{(\alpha)}((a_i)_{1\leq i \leq p}, (b_j)_{1 \leq j \leq q}; x,y) = \sum_{n=0}^{\infty}\sum_{|\tau| = n}\frac{(a_1)_{\tau}\cdots(a_p)_{\tau}}{(b_1)_{\tau}\cdots(b_q)_{\tau}}
\frac{C_{\tau}^{(\alpha)}(x)C_{\tau}^{(\alpha)}(y)}{C_{\tau}^{(\alpha)}(1)|\tau|!},
\end{gather*}
where $(1)= (1,\dots,1)$ and for a partition $\tau$
\begin{gather*}
(a)_{\tau}^{(\alpha)} = \prod_{i=1}^m \left(a - \frac{1}{\alpha}(i-1)\right)_{\tau_j}= \prod_{i=1}^m \frac{\Gamma(a - (i-1)/\alpha + \tau_j)}{\Gamma(a - (i-1)/\alpha)}
\end{gather*}
is the generalized Pochhammer symbol and $\Gamma$ is the Gamma function. We assume of course that the above expression makes sense for all the coef\/f\/icients
$(a_i)_{1 \leq i \leq p}$, $(b_j)_{1 \leq j \leq q}$. An interesting feature of the hypergeometric series of Jack parameter $\alpha = 1$ is that they are expressed through univariate  hypergeometric functions ${}_p\mathscr{F}_q$~\cite{Leb} as follows~\cite{Gro}
\begin{gather*}
 {}_pF_q^{(1)}((m+\mu_i)_{1\leq i \leq p}, (m+\phi_j)_{1\leq j\leq q};x, y) =\pi^{\frac{m(m-1)}{2}(p-q-1 + 1/\alpha)} \prod_{i=1}^m (m - 1)!
\\\qquad{}\times \prod_{i=1}^p\frac{(\Gamma(\mu_i+1))^m}{\Gamma(m+\mu_i)} \prod_{j=1}^q\frac{\Gamma(m+\phi_j)}{(\Gamma(\phi_j+1))^m}
\frac{\det\big[{}_p\mathscr{F}_q((\mu_i+1)_{1\leq i \leq p}, (1+\phi_j)_{1\leq j\leq q};x_ly_r )\big]_{l,r=1}^m}{V(x)V(y)}
\end{gather*}
for all $\mu_i,\phi_j > -1$, where $V$ stands for the Vandermonde polynomial and empty product are equal $1$. For instance, when $p=q=0$, we have \begin{align*}
{}_0F_0^{(1)}(x, y) =\pi^{-\frac{m(m-1)}{2}(-1+1/\alpha)} \frac{\det(e^{x_iy_j})_{i,j=1}^m}{V(x)V(y)},
\end{align*}
and for $p=0$, $q=1$, we similarly get
\begin{gather*}
 {}_0F_1^{(1)} (m+\phi; x, y) =\pi^{\frac{m(m-1)}{2\alpha}} \prod_{i=1}^m (m - 1)! \frac{\Gamma(m+\phi)}{(\Gamma(\phi+1))^m}\frac{\det\big[{}_0\mathscr{F}_1(1+\phi; x_iy_j)\big]_{i,j=1}^m}{V(x)V(y)}
\end{gather*}
for $\phi > -1$.

\section[Brownian motions in Weyl chambers: types $A$ and $B$ revisited]{Brownian motions in Weyl chambers: types $\boldsymbol{A}$ and $\boldsymbol{B}$ revisited}

Brownian motions in Weyl chambers were deeply studied in \cite{Gra} and they are shown to be $h$-transforms in Doob's sense of $m$ independent real BMs killed when they f\/irst hit $\partial C$. They are then interpreted as $m$ independent particles constrained to stay in the Weyl chamber $C$. The proofs of those properties use probabilistic arguments. Here, we show how, for both types~$A$ and~$B$, these properties follow from the above determinantal representation of multivariate hypergeometric functions of two arguments.
Let us f\/irst recall that the $A$-type root system is def\/ined by
\begin{gather*}
A_{m-1} = \{\pm (e_i - e_j), \  1 \leq i < j \leq m\}, \end{gather*} with positive and simple systems given by
\begin{gather*}
R_+ = \{e_i - e_j, \ 1\leq i < j \leq m \}, \qquad S = \{e_i - e_{i+1}, \ 1\leq i  \leq m-1 \},
\end{gather*}
where $(e_i)_{1 \leq i \leq m}$ is the standard basis of $\re^m$, $W$ being the permutation group. In this case, $V = \re^m$, the span of $R$ is the hyperplane of $\re^m$ consisting of vectors whose coordinates sum to zero and $C = \{x \in \re^m, x_1 > \dots > x_m\}$. Besides, there is only one orbit so that $k(\alpha) := k_1 \geq 0$ and the function $h$ (product of the positive roots) is given by the Vandermonde polynomial.

Next,  it was shown in \cite{Bak} that $D_k^W$ is given up to the constant $1/|W|$ by ${}_0F_0^{(1/k_1)}(x,y)$ (we use another normalization than the one used in \cite{Bak} so that the factor $\sqrt{2}$ is removed). Thus, when $k_1=1$, (\ref{sg}) writes for $R=A_{m-1}$
\begin{gather}\label{BMAWC}
p_t^{1}(x,y) = \frac{V(y)}{V(x)}\det(N_t(y_j-x_i))_{i,j=1}^m,
\end{gather}
where $N_t$ is the heat kernel, namely
\begin{gather*}
N_t(v) = \frac{1}{\sqrt{2\pi t}}e^{-v^2/2t}.
\end{gather*}
Thus, the $V$-transform property is easily seen. A similar result holds for $R = B_m$. This root system has the following data
\begin{alignat*}{3}
& R  =  \{\pm e_i, \ \pm e_i \pm e_j,\  1 \leq i < j \leq m\},\qquad
   & &  R_+  =  \{e_i, \ 1\leq i \leq m, \  e_i \pm e_j,\  1 \leq i < j \leq m\},&\\
& S  =  \{e_i - e_{i+1}, \  1 \leq i \leq m, \  e_m \},\qquad & & C = \{y \in \re^m, \  y_1 > y_2 > \cdots > y_m > 0\}.
\end{alignat*}
$W$ is generated by transpositions and sign changes $(x_i \mapsto -x_i)$ and there are two orbits so that $k=(k_0,k_1)$ thereby $\gamma = mk_0 + m(m-1)k_1$. The generalized Bessel function\footnote{There is an erroneous sign in one of the arguments in \cite{Bak}. Moreover, to recover this expression in the $B_m$ case from that given in \cite{Bak}, one should make substitutions $a = k_0 - 1/2$, $k_1 = 1/\alpha$, $q = 1 +(m-1)k_1$.} is given by~\cite[p.~214]{Bak}
\begin{gather*}
\frac{1}{|W|}D_k^W(x,y) = {}_0F_1^{(1/k_1)}\left(k_0 + (m-1)k_1 + \frac{1}{2}, \frac{x^2}{2},\frac{y^2}{2}\right).
\end{gather*}
It follows that
\begin{align*}
p_t^{k_0,1}(x,y) & = C(m)\frac{h(y)}{h(x)}\frac{e^{-(|x|^2+|y|^2)/2t}}{t^{m/2}}\prod_{i,j=1}^m\left(\frac{x_iy_j}{t}\right)\det\left[\be\left(k_0 + \frac{1}{2}, \frac{(x_iy_j)^2}{4t^2}\right)\right]
_{i,j=1}^n,\end{align*}
where $h(y) = V(y^2)\prod\limits_{i=1}^ny_i$ and $C(m)$ is a constant depending on $m$. Taking $k_0 = 1$ and using the identity~\cite{Leb}
\begin{gather*}
\be\left(\frac{3}{2},z\right) =  \frac{C}{2\sqrt{z}}\sinh\big(2\sqrt{z}\big),
\end{gather*}
for some constant $C$, one gets for the Brownian motion in the $B$-type Weyl chamber
\begin{align}\label{BMBWC}
p_t^{1,1}(x,y)  =  \frac{h(y)}{h(x)}\det\left[N_t(y_j - x_i) - N_t(y_j+x_i)\right]_{i,j=1}^m.
\end{align}
The $h$-transform property is then obvious.

\section[Generalized Bessel function of type $D$]{Generalized Bessel function of type $\boldsymbol{D}$}

The root system of type $D$ is def\/ined by \cite[p.~42]{Hum}
\begin{gather*}
R  = \{\pm e_i \pm e_j, \ 1 \leq i < j \leq m\},\qquad R_+ = \{e_i \pm e_j, \ 1 \leq i < j \leq m\},
\end{gather*}
and there is one orbit so that $k(\alpha) = k_1$ thereby $\gamma = m(m-1)k_1$. The Weyl chamber is given by:
\begin{gather*}
C = \{x \in \re^m,\  x_1 > x_2 > \dots > |x_m|\} = C_1 \cup s_m C_1,
\end{gather*}
where $C_1$ is the Weyl chamber of type $B$ and $s_m$ stands for the ref\/lection with respect to the vector $e_m$ acting by sign change on the variable $x_m$.
Now, Grabiner's result reads in this case
\begin{align*}
p_t^1(x,y) &= \frac{V(y^2)}{V(x^2)}\frac{\det[N_t(y_i - x_j) - N_t(y_i + x_j)]_{i,j=1}^m + \det[N_t(y_i - x_j) + N_t(y_i + x_j)]_{i,j=1}^m}{2}
\\ & = \frac{C_m}{t^{\gamma + m/2}}e^{-(|x|^2+|y|^2)/2t}\frac{\det\left[\sinh(x_iy_j/t)\right]_{i,j=1}^m + \det\left[\cosh(x_iy_j/t)\right]_{i,j=1}^m}{V(x^2/4t^2)V(y^2)}V^2(y^2),
\end{align*}
where $\gamma = m(m-1)$. With the help of the determinantal representations  and of~\cite{Leb}
\begin{gather*}
\be\left(\frac{3}{2},z\right) =  \frac{C}{\sqrt{z}}\sinh\big(2\sqrt{z}\big), \qquad \be\left(\frac 12,  z\right) = \cosh\big(2\sqrt{z}\big),
\end{gather*}
one gets
\begin{gather*}
p_t^1(x,y) =  \frac{e^{-(|x|^2+|y|^2)/2t}}{c_k t^{\gamma + m/2}}\\
\phantom{p_t^1(x,y) =}{} \times
\left[\prod_{i=1}^m \left(\frac{x_iy_i}{2t}\right){}_0F_1^{(1)}\left(m+\frac{1}{2}, \frac{x^2}{2t},\frac{y^2}{2t}\right)+ {}_0F_1^{(1)}\left(m-\frac{1}{2}, \frac{x^2}{2t},\frac{y^2}{2t}\right)\right]
V^2(y^2).
\end{gather*}
With regard to (\ref{sg}) and setting $q := 1+(m-1)k_1$, it is natural to prove the claim of Theorem~\ref{theorem1}.

\begin{proof} It uses the so-called {\it shift principle} that we brief\/ly outline~\cite{Dunkl}. Let $E$ be a conjugacy class of roots of $R$ under the action $W$. Let $k_E$ be the value of the multiplicity function on this class. Then, the generalized Bessel function associated with the multiplicity function $k'$ def\/ined by
\begin{gather*}
k'(\alpha) = \left\{\begin{array}{ll}
k(\alpha) + 1 & \textrm{if} \ \ \alpha \in E, \\
k(\alpha) & \textrm{otherwise},
\end{array}\right.
\end{gather*}
is given by
\begin{gather*}
\frac{1}{|W|}D_{k'}^W(x,y) = C\sum_{w \in W}\xi_E(w)\frac{D_k(x,wy)}{p_E(x)p_E(y)},\qquad p_E(x) := \prod_{\alpha \in R_+ \cap E}\langle \alpha,x\rangle,
\end{gather*}
where $\xi_{E}(w)$ is def\/ined by $p_E(wx) = \xi_E(w)p_E(x)$ and $C$ is some constant. Since we will deal with both the $B$ and $D$-types, it is convenient to add a superscript $B$ or $D$ to each corresponding item. Therefore, $W^B$, $W^D$ denote the Weyl groups associated with the root systems of types~$B$,~$D$ respectively and $k^B$, $k^D$ denote the corresponding multiplicity functions. Recall that~\cite{Hum} $W^B$ is the semi-direct product of $S_m$ and $(\mathbb{Z}/2\mathbb{Z})^m$ (sign changes) while $W^D$ is the semi-direct product of $S_m$ and $(\mathbb{Z}/2\mathbb{Z})^{m-1}$ (even sign changes). Let $E := \{e_i,\ 1 \leq i \leq m\}$, then $p_E(x)$ is invariant under permutations and even sign changes and skew-invariant under odd sign changes
(note that $D_k^{W^D}$ is not $W^B$-invariant). It follows that $\xi_E(w) = 1$ for $w \in W^D$ while $\xi_E(w) = -1$ for $w \in W^B\setminus W^D$. Now, recall that $k^D$ takes only one value while $k^B$ takes two values. Since the Dunkl Laplacian $\Delta_k^B$ of $B$ reduces to $\Delta_k^D$ when the value of the multiplicity function $k_0$ on~$E$ is zero~\cite{Ros1}, the (non-symmetric) Dunkl kernel $D_k^D$ associated with $R=D_m$ is given by the Dunkl kernel $D_k^B$ associated with $R=B_m$ when $k_E^B = 0$.
In fact, this is true since $D_k$ is the solution of a spectral problem which is independent of $W$.
As a result,
\begin{align*}
\frac{1}{|W^D|}D_k^{W^D}(x,y) &= \frac{1}{|W|^D}\sum_{w \in W^D}D_k^{D}(x,wy)\\&
= \frac{1}{2|W|^D}\sum_{w \in W^B}(1+ \xi_E(w))D_k^{B, k_0 = 0}(x,wy)\\&
= \frac{|W^B|}{2|W^D|}\left[\frac{1}{|W^B|}D_{k}^{W^B (k_0 = 0)} + \frac{p_E(x)p_E(y)}{C|W^B|}D_{k'}^{W^B (k_0=1)}\right],
\end{align*}
where we used the shift principle to derive the last line. Keeping in mind that
\begin{gather*}
\frac{1}{|W^B|}D_k^{W^B}(x,y) = {}_0F_1^{(1/k_1)}\left(k_0 + (m-1)k_1 + \frac{1}{2}, \frac{x^2}{2t},\frac{y^2}{2t}\right)
\end{gather*}
and using that $|W^B| = 2|W^D|$ \cite[p.~44]{Hum}  we are done with $C= 2^m$.
\end{proof}

\section{Concluding remarks}
As the reader can see, $D_k^{W^D}$ is not equal to $D_k^{W^B}$ specialized with $k_0 = 0$ and this shows that the generalized Bessel function is intimately related to $W$. Moreover the symmetrical of $D_k^{W^D}$ with respect to $s_m$ gives $D_k^{W^B}$ in the special case $k_0=0$ which ref\/lects the fact that the Weyl chamber of type $D$ is the union of the one of type $B$ and its symmetrical with respect to $s_m$.

\pdfbookmark[1]{References}{ref}
\LastPageEnding

\end{document}